\documentstyle[11pt]{article}
\renewcommand{\baselinestretch}{1.3}
\newenvironment {qth}{\\[2ex] {\bf Theorem:} \it}{{ }\\}

\def\Cox{\hfill \Box}

\def\iid{\mbox{i.i.d.}}

\def\E{{\bf{E}}}
\def\P{{\bf{P}}}

\def\R{{\bf{R}}}
\def\Z{{\bf{Z}}}

\def\|{\, | \, }

\begin{document}
 
\setcounter{equation}{0}

\begin{titlepage}
\begin{center}
{\large \bf ON WHICH GRAPHS ARE ALL RANDOM WALKS IN RANDOM
ENVIRONMENTS TRANSIENT?} \\
\end{center}
\vspace{5ex}
\begin{flushright}
Robin Pemantle \footnote{Research supported in part by 
National Science Foundation grant \# DMS 9300191, by a Sloan Foundation
Fellowship, and by a Presidential Faculty Fellowship}$^,$\footnote{Department 
of Mathematics, University of Wisconsin-Madison, Van Vleck Hall, 480 Lincoln
Drive, Madison, WI 53706}
  ~\\
Yuval Peres \footnote{Department of
Statistics, University of California, Berkeley, CA 94720} 
\end{flushright}

\vfill

{\bf SUMMARY:} \break
An infinite graph $G$ has the property that a random walk in 
random environment on $G$ defined by i.i.d.\ resistances
with any common distribution is almost surely transient, 
if and only if for some $p < 1$, simple random walk is transient 
on a percolation cluster of $G$ under bond percolation with
parameter $p$.
\vfill



\end{titlepage}

Simple random walks on infinite graphs have been widely studied
with particular attention to the dichotomy between recurrence
and transience.  In the last two decades there has been interest
in random walks in random environments (RWRE's), which are mixtures
of Markov chains in which the transition probabilities are
picked according to some prior distribution.  For reversible
Markov chains, the transition probabilities may be described
by associating resistances $R(vw)$ to each edge $(vw)$, 
the transition probabilities out of any vertex being 
proportional to the inverses of the resistances on the edges 
incident to that vertex: 
\begin{equation} \label{eq hi its late}
p (v,w) = { R(vw)^{-1} \over \sum_z R(vz)^{-1}}
\end{equation}
where the sum is over all vertices $z$ adjacent to $v$.
One natural model\footnote{see Remark~1 for another RWRE model
relevant to the present discussion},
considered in Grimmett and Kesten (1984),
involves taking the resistances as $\iid$ positive random variables.  
A result of Adams and Lyons (1991) implies 
that if a tree $\Gamma$ has positive Hausdorff dimension,
then for an arbitrary distribution of the resistances 
the resulting network is almost surely transient.  
(See Lyons (1990) where the term {\em branching number} is
used for the exponential of the Hausdorff dimension.) 
The surprising aspect of this is that no conditions 
are imposed on the tail of the distribution.  
This note proves a converse of this implication, 
but on a general graph.  Specifically, the implication
$(a) \Rightarrow (b)$ in part~$(ii)$ of the theorem below
shows that almost sure transience for an arbitrary 
resistance distribution on a tree $\Gamma$ implies $p_c < 1$ 
for IID bond percolation.  By results in Lyons (1990), this
implies that $\Gamma$ has positive dimension.
Part~$(i)$ is immediate from the well known
connections between random walks and electrical networks
(c.f.\ Doyle and Snell (1984)) and is included for completeness.
\begin{qth}
Let $G$ be a locally finite connected graph.  Suppose that a positive
random variable $R(e)$ is attached to every edge $e$ of $G$ so that
the collection $\{ R(e) \}$ is $\iid$.  Let $\mu$ denote their 
common distribution.  Consider the RWRE given by~(\ref{eq hi its late}).  
\begin{quote}
$(i)$  If $\int x \, d\mu < \infty$ and simple random walk on $G$ is
transient, then the resistor network $\{ R(e) \}$ has finite net resistance
almost surely, thus defining almost surely a transient random walk.

$(ii)$  The following conditions on $G$ are equivalent:
$(a)$ for any distribution $\mu$ on $(0,\infty)$ the resistor network
$\{ R(e) \}$ is almost surely transient; $(b)$ there exists some
$p < 1$ such that independent bond percolation on $G$ with parameter~$p$ 
contains with positive probability a cluster on which 
simple random walk is transient.
\end{quote}
~~
\end{qth}

{\sc Proof:}  We rely on two facts from the theory relating 
reversible random walks to the corresponding resistor networks.  
Firstly, the random walk is transient if and only if the network
has finite resistance from some (hence every) vertex to infinity, 
i.e., there is some nonzero flow $\{ F(e) \}$, satisfying Kirchoff's 
laws except at one source vertex, such that the energy
$$\sum_e F(e)^2 R(e)$$
is finite.  Secondly, Rayleigh's principle asserts that
the resistance of a network is monotone in each of the 
individual resistances.  Both facts may be found in Doyle
and Snell (1984).  The proof of~$(i)$ follows immediately:
$G$ is transient for simple random walk, so some flow
$\{ F(e) \}$ on $G$ has finite energy 
$$ \sum_e F(e)^2$$
with respect to unit resistances; the expected energy of this flow 
in the network of random resistances is
$$\E \sum_e R(e) F(e)^2 = \int x \, d\mu \sum_e F(e)^2 < \infty $$
so the resistance is almost surely finite and hence 
the network is almost surely transient.

$(ii)$  The direction $(b) \Rightarrow (a)$ is easy: assume $(b)$
for some parameter $p$ and consider the random subgraph of $G$
consisting of those edges $e$ for which $R(e) \leq Q_p$, where
$1 > \mu (0,Q_p) \geq p$.  Resistances have the following linearity:
multiplying every resistor by $C$ multiplies the resistance of the 
network by $C$.  Thus our assumption~$(b)$ implies there is a 
positive probability that some component of this subgraph 
will have finite resistances when all edges have resistance $Q_p$.  
By Rayleigh's law, decreasing the resistances on edges in
this component from $Q_p$
to $R(e)$ and decreasing all other resistances from $\infty$
to $R(e)$ gives a network with finite resistance, and therefore
random walk on $G$ with resistances $\{ R(e) \}$ is transient
with positive probability.  Since this is a 
tail event for the countable collection $\{ R(e) \}$, 
the probability is one.

Now we verify that $(a) \Rightarrow (b)$.  Assume that for every
$p < 1$, simple random walk on $p$-percolation clusters is almost 
surely recurrent.  Let $p_k \uparrow 1$ and construct a probability 
measure $\mu$ on $(0,\infty)$ as follows.  Fix a vertex $v \in G$
and set $\gamma_1 = 1$.  Define a measure $\mu_1$ on $\R \cup \{ \infty \}$
by $\mu_1 (\{\gamma_1\}) = p_1$ and $\mu_1 (\{\infty \} ) = 1 - p_1$.  
The RWRE on $G$ with resistance distribution $\mu_1$ 
is just simple random walk on the $p_1$-percolation clusters
and is thus recurrent.  Therefore there exists $N_1$ such that
when choosing a random environment according to $\mu_1$, the probability
is at least $1/2$ that either all edges adjacent to $v$ have infinite
resistance (i.e. $v$ is isolated) or else the resulting RWRE started from $v$ 
returns to $v$ in the first $N_1$ steps.  Let $D_1$ be the maximal degree of 
any vertex in $G$ within distance $N_1$ of $v$.

For the induction step, let $k \geq 2$ and assume that $\gamma_j ,
\mu_j , N_j$ and $D_j$ have been defined for $1 \leq j < k$ and
$\mu_j (\infty) = 1 - p_j$.  Pick $\gamma_k > 2 N_{k-1}
D_{k-1} \gamma_{k-1}$.  Define $\mu_k$ by moving mass $p_k - p_{k-1}$ 
from infinity to $\gamma_k$, i.e. 
$$\mu_k = \mu_{k-1} + (p_k - p_{k-1}) (\delta_{\gamma_k} - \delta_\infty) .$$
The RWRE with resistance distribution $\mu_k$ is recurrent since
it yields an environment whose resistances on the
$p_k$-percolation cluster are bounded between 1 and $\gamma_k$
[use the aforementioned linearity together
with Rayleigh's law].  Thus we may choose an integer $N_k$, 
such that if a random environment is chosen
according to $\mu_k$, the probability is at most
$2^{-k}$ that $v$ is not isolated and the resulting RWRE started
from $v$ fails to return to $v$ within the first $N_k$ steps.  
Finally, define $D_k$ to be the maximal degree of any vertex of $G$ 
within distance $N_k$ of $v$.  

Finally, set $\mu (\{ \gamma_k \}) = p_k - p_{k-1}$ for all $k \geq 2$,
with $\mu (\{ \gamma_1 \}) = p_1$, so that $\mu$ is the weak limit
of the $\mu_k$.  We construct a $\mu$-RWRE $\{ S_n \}$ 
on the same probability space as a sequence
$\left ( \{ S_n^{(k)} \} : k = 1,2,3,\ldots \right )$ of
$\mu_k$-RWRE's as follows.  Let the resistances $\{ R(e) \}$ 
be i.i.d.\ with common distribution $\mu$ and define
$$R^{(k)} (e) = \left \{ \begin{array}{ll} R(e) &
   \mbox{if } R(e) \leq \gamma_k \\ \infty & 
   \mbox{if } R(e) > \gamma_k \end{array} \right. \; . $$
Observe that this makes $\{ R^{(k)} (e) \}$ i.i.d.\ with
distribution $\mu_k$.
Let $\{ S_n \}$ be a random walk starting from $v$
with transition probabilities determined by the resistances 
$\{ R(e) \}$.  Note that if $R(wz) \leq \gamma_k$ then
the transition probability from $w$ to $z$ with resistances
$\{ R^{(k)} (e) \}$ is greater than or equal to the
the transition probability from $w$ to $z$ with resistances
$\{ R(e) \}$.  Thus for each $k$, we may define a random sequence 
$\{ S_n^{(k)} \}$ so as to have transition probabilities 
determined by $\{ R^{(k)} \}$ and so that $S_n^{(k)} = S_n$ for all
$n \leq T_k$, where $T_k$ is the least time $t$ for which
$R(S_t S_{t+1}) > \gamma_k$.  

Define events 
\begin{eqnarray*}
A_k & = &  \{ R(vw) > \gamma_k \mbox{ for all } w \mbox{ neighboring } v 
   \} ~\subseteq~ \{ T_k = 0 \} \\
B_k & = & \{ R(S_j,S_{j+1}) \geq \gamma_{k+1} \mbox{ for some } 
   0 \leq j < N_k \} \setminus A_k ~=~ \cup \{ 0 \leq T_k < N_k \}
   \setminus A_k \\
C_k & = & \{ S_n^{(k)} \neq v \mbox{ for all } 1 \leq n \leq N_k \} .
\end{eqnarray*}
The event $G_k = \{ S_n \neq v \mbox{ for all } 1 \leq n \leq N_k \}$
is contained
in $A_k \cup B_k \cup C_k$, since on $G_k \setminus C_k$ the
time $T_k$ is less than $N_k$.  For $1 \leq j < N_k$, the probability of 
the event $\{ T_k = j \}$ is at most $D_k \gamma_k / \gamma_{k+1}$,
since the sum of $R(e)^{-1}$ over edges incident to $S_j$ 
with resistance greater than $\gamma_k$ is at most $D_k / \gamma_{k+1}$
and there is at least one edge incident to $S_j$ 
with $R(e)^{-1} \geq \gamma_k^{-1}$.  Similarly, $\P (\{ T_k = 0 \} 
\setminus A_k) \leq D_k \gamma_k / \gamma_{k+1}$.  
Thus $\P (B_k) \leq N_k D_k \gamma_k / \gamma_{k+1} \leq 2^{-k}$
by construction of $\gamma_{k+1}$.  But $N_k$ is defined so
that $\P (C_k \setminus A_k) \leq 2^{-k}$ and clearly $\P (A_k) 
\leq 1 - p_k$.  Therefore
$$\P (G_k) \leq (1 - p_k) + 2^{1-k} \rightarrow 0$$
as $k \rightarrow \infty$ and the $\mu$-RWRE is recurrent.   $\Cox$

{\bf Remarks:} \\

\noindent{1. } The recent paper of Grimmett, Kesten and Zhang (1991)
shows that for Euclidean lattices $\Z^d$ ($d \geq 3$), random walk
on supercritical percolation clusters is almost surely transient,
so the theorem is applicable, and all RWRE's with i.i.d.\ positive
resistances are transient.  Note that the RWRE models discussed 
in Durrett (1986) and Sunyach (1987) may be recurrent even when 
$d \geq 3$.  Resistances in Durrett's model are stationary under
$\Z^d$-shifts, so his model is close to the present one.

\noindent{2. } By Lyons (1990), the critical probability for bond
percolation on a tree $\Gamma$ is $p_c (\Gamma) = e^{-dim \Gamma}$ 
and if $dim (\Gamma) > 0$ and $p > p_c$, then
simple random walk on any infinite cluster is transient
(since the cluster must have positive dimension).  For a tree,
property $(b)$ is thus equivalent to $dim (\Gamma) > 0$.
The Max Flow Min Cut Theorem (as applied in Lyons 1990) allows us 
to formulate this purely in terms of flows: 
\begin{quote}
On a tree, property $(a)$ is equivalent to the existence of
a flow $\{ F(e) \}$ satisfying Kirchoff's laws (except at the root) 
for unit resistors, and having the exponential decay property: 
$$F(e) \leq C \rho^{\mbox{dist}(e)}$$
where $\mbox{dist}(e)$ is the distance from $e$ to the root of $\Gamma$.
\end{quote}

Acknowledgement: Thanks to Harry Kesten for showing us a 
surprising example which led to this note.

\renewcommand{\baselinestretch}{1.0}\large


\normalsize
\begin{thebibliography}{YMN}

\bibitem{AL}
Adams, S. and Lyons, R. (1991).  Amenability, Khazdan's property and
percolation for trees, groups and equivalence relations.  {\em Israel
J. Math.} {\bf 75} 341 - 370.

\bibitem{DS}
Doyle, P. and Snell, J. L. (1984).  Random walks and electrical networks.
Mathematical Association of America: Washington.

\bibitem{Du}
Durrett, R. (1986).  Multidimensional random walks in random environments
with subclassical limiting behavior.  {\em Comm. Math. Phys.} {\bf 104}
87 - 102.

\bibitem{GK}
Grimmett, G. and Kesten, H. (1984).  Random electrical networks on complete
graphs.  {\em J. London Math. Soc.} {\bf 2}  171 - 192.

\bibitem{GKZ}
Grimmett, G., Kesten, H. and Zhang, Y. (1991).  Random walk on the
infinite cluster of the percolation model.  {\em Preprint.}

\bibitem{Ly}
Lyons, R. (1990).  Random walks and percolation on trees.  {\em Ann.
Probab.} {\bf 18} 931 - 958.

\bibitem{Su}
Sunyach, C. (1987).  Sur la transience et la r\'ecurrence des marches
al\'eatoires en milieu al\'eatoire.  {\em Ann. Inst. H. Poin.} {\bf 23}
613 - 626.

\end{thebibliography}
\end{document}